\documentclass[english]{article}
\usepackage[T1]{fontenc}
\usepackage[latin9]{inputenc}
\usepackage{amsmath}
\usepackage{graphicx}
\usepackage{amssymb}

\newcommand{\lyxline}[1][1pt]{%
  \par\noindent%
  \rule[.5ex]{\linewidth}{#1}\par}

\usepackage{babel}

\begin{document}

\title{Yet Another Riemann Hypothesis}

\author{Linas Vepstas <linasvepstas@gmail.com>}

\date{12 October 2004; revised 31 December 2010}
\maketitle
\begin{abstract}
This short note presents a peculiar generalization of the Riemann
hypothesis, as the action of the permutation group on the elements
of continued fractions. The problem is difficult to attack through
traditional analytic techniques, and thus this note focuses on providing
a numerical survey. These results indicate a broad class of previously
unexamined functions may obey the Riemann hypothesis in general, and
even share the non-trivial zeros in particular.
\end{abstract}

\section{Introduction}

The Riemann zeta can be expressed as an integral in the following
form\cite{Edw74}:\[
\zeta(s)=\frac{s}{s-1}-s\int_{0}^{1}h(x)x^{s-1}dx\]
where \[
h(x)=\frac{1}{x}-\left\lfloor \frac{1}{x}\right\rfloor \]
is sometimes called the Gauss map. Here, $\left\lfloor y\right\rfloor $
is the largest integer less than or equal to $y$. The Gauss map lops
of a digit in the continued fraction expansion of $x$. If one writes
the continued fraction expansion\cite{Khin35} of $0<x\le1$ as\[
x=[a_{1},a_{2},a_{3},\cdots]=\frac{1}{a_{1}+\frac{1}{a_{2}+\frac{1}{a_{3}+\ldots}}}\]
then $h(x)=[a_{2},a_{3},\cdots]$. It is a (reverse) shift operator
on the continued fraction expansion. The study of continued fractions
is interesting because of their proximity to various problems of fractals
and dynamical systems, including ties to Farey fractions, the modular
group via the Cantor set, and the Minkowski question mark function.
More immediately, the transfer (Ruelle-Frobenius-Perron) operator
of the Gauss map is the Gauss-Kuzmin-Wirsing operator, which has a
number of interesting properties of its own.

In this note, we observe that the Gauss map can also be thought of
as one particular element of the permutation group acting on an infinite
dimensional representation of the real numbers. Thus, we are lead
to contemplate other elements of the permutation group. Specifically,
consider the permutation operator \begin{equation}
S_{p,q}:[a_{1},\cdots,a_{p},\cdots,a_{q},\cdots]\mapsto[a_{1},\cdots,a_{q},\cdots,a_{p},\cdots]\label{eq:Permutation Elts}\end{equation}
which exchanges the $p$'th and $q$'th digits of the continued fraction
expansion of $x$. Note that $S_{p,q}$ maps the unit interval onto
the unit interval, and is discontinuous at a countably infinite number
of points. The generalized Riemann hypothesis that we choose to explore
concerns the zeros of the integral\begin{equation}
\zeta_{p,q}(s)=\frac{s}{s-1}-s\int_{0}^{1}S_{p,q}(x)x^{s-1}dx\label{eq:zeta-permute}\end{equation}
That is, is it possible that the (non-trivial) zeros of $\zeta_{p,q}(s)$
all lie on the $\Re s=\frac{1}{2}$ axis? This text is devoted to
the exploration of this possibility. 

This integral can be evaluated, but only with considerable difficulty.
The function $S_{p,q}\left(x\right)$ is a piece-wise combination
of unimodular Mobius functions, that is, a piece-wise collection of
\[
\frac{ax+b}{cx+d}\]
with integer values of $a$, $b$, $c$ and $d$, such that $ad-bc=\pm1$.
This structure follows from the general theory of continued fractions\cite{Khin35}.
Section \ref{sec:Structural overview} below provides a set of graphs
visualizing $S_{p,q}(x)$, showing their piece-wise nature. The functions
are clearly {}``fractal'' or self-similar; a general discussion
of the self-similarity is given in \cite{Vep-mink2008}. Because of
the piece-wise structure, the integral can be evaluated analytically.
The result is a messy set of nested sums that don't seem to provide
any particular insight. One might expect that these sums would simplify,
or at least, speed up numeric evaluation, but they don't even seem
very useful for that. Details are given in section \ref{sec:Sums-and-Integrals}.

It is easiest, at first, to perform the integration numerically. Treating
the integral as a very simple Newtonian integration sum with equally
spaced abscissas, one gets:\begin{equation}
\zeta_{N;p,q}\left(s\right)=\frac{s}{s-1}-\frac{s}{N}\;\sum_{n=0}^{N-1}S_{p,q}\left(\frac{2n+1}{2N}\right)\;\left(\frac{2n+1}{2N}\right)^{s-1}\label{eq:zeta-pq-sum}\end{equation}
This summation converges very slowly and very noisily as $N\to\infty$
but a numerical evaluation for large $N$ can give some basic insight.
First and foremost, one discovers that, at least qualitatively, the
sums $\zeta_{N;p,q}\left(s\right)$ behave very much like the corresponding
sum for the Riemann zeta. This provides the simplest evidence in support
of the hypothesis. A graphical visualization is provided in section
\ref{sec:Structural overview}.

\section{Continued Fractions; Permutation Groups}

Continued fractions provide a representation of the real numbers in
the infinite Cartesian product space $\mathbb{M}^{\omega}=\mathbb{M\times M\times}...$
where $\mathbb{M}=\mathbb{N}\cup\{\infty\}$ and $\mathbb{N}=\mathbb{Z}^{+}$
is the set of positive integers. The continued fraction brackets are
a function that maps this product space to the real numbers, \emph{i.e.}
$[\,]:\mathbb{N}^{\omega}\mathbb{\rightarrow R}$. The mapping is,
strictly-speaking, surjective, yet, in a certain sense, is {}``almost
everywhere'' bijective. That is, every non-rational real number corresponds
to a unique continued fraction; only the rationals have multiple expansions. 

The above definition of a continued fraction differs slightly from
the conventional norm. Usually, continued fractions are considered
to be of finite length for the rationals, and infinite-length for
the non-rational reals. This can create confusion when discussing
$S_{p,q}(x)$ when either $p$ or $q$ are larger than the length
of the continued fraction. The trick of adjoining infinity to the
naturals is the most straight-forward way of avoiding this problem.
The cost of this trick is that a {}``very large part'' of the Cartesian
product maps to the rationals. Thus, for example, the continued fraction
expansion for 0 is $0=[\infty,a_{2},a_{3},...]$ for any arbitrary
positive integers $a_{k}$. All rationals suffer in the same way:
a finite continued fraction is just one where $\infty$ appears somewhere
in the expansion, i.e. $[a_{1},a_{2},...,a_{N}]=[a_{1},a_{2},...,a_{N},\infty,a_{N+2},...]$.
This can be papered over by considering the reals as the quotient
space of sequences modulo the kernel of $[\,]$. That is, define a
quotient space $\mathbb{S}=\mathbb{M}^{\omega}/\ker[\,]$, so that
the operator $[\,]:\mathbb{S}\to\mathbb{R}$ is 1-1 and onto in this
quotient space. Since this map is bijective, it is invertible, and
thus, any function on the reals $\mathbb{R}$ and be expressed as
an equivalent function on $\mathbb{S}$, and vice-versa: functions
on the sequence space $\mathbb{M}^{\omega}$ correspond to functions
on the reals. 

Elements of the (infinite-dimensional) permutation group acting on
$\mathbb{S}$ are given by equation \ref{eq:Permutation Elts}. The
$S_{p,q}(x)$ are discontinuous, but only on a countable set of points.
That is, for fixed $p,q,$ the function $S_{p,q}(x)$ is discontinuous
on the proverbial ''set of measure zero'', and so one can easily
form well-defined integrals using ordinary techniques.

\section{Structural overview\label{sec:Structural overview}}

A graphical presentation of the functions $S_{1,2}\left(x\right)$,
$S_{1,3}\left(x\right)$, $S_{1,4}\left(x\right)$ and $S_{2,3}\left(x\right)$
are shown in figure \ref{fig:Continued-Fraction-Permutations}. The
fractal, self-similar nature of these functions is readily apparent.

\begin{figure}
\caption{Continued Fraction Permutations\label{fig:Continued-Fraction-Permutations}}

\includegraphics[width=0.5\textwidth]{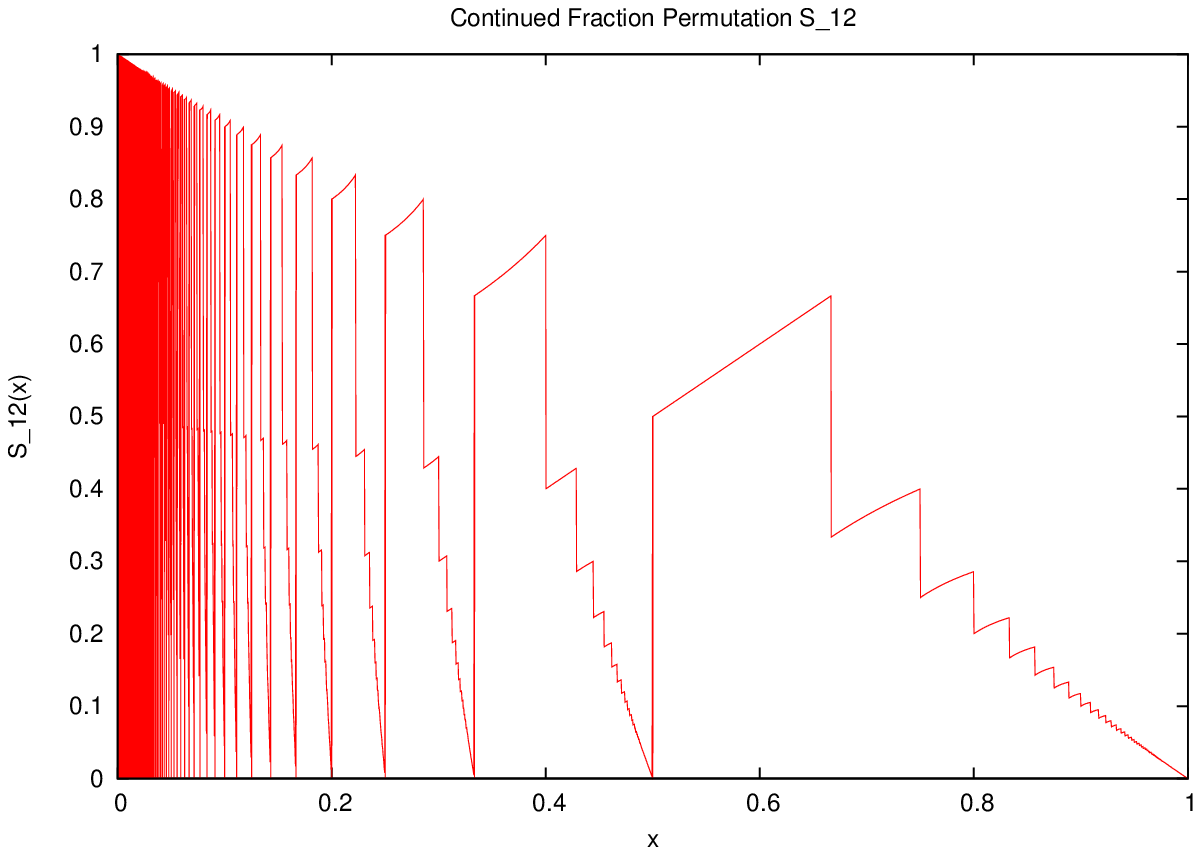}\includegraphics[width=0.5\textwidth]{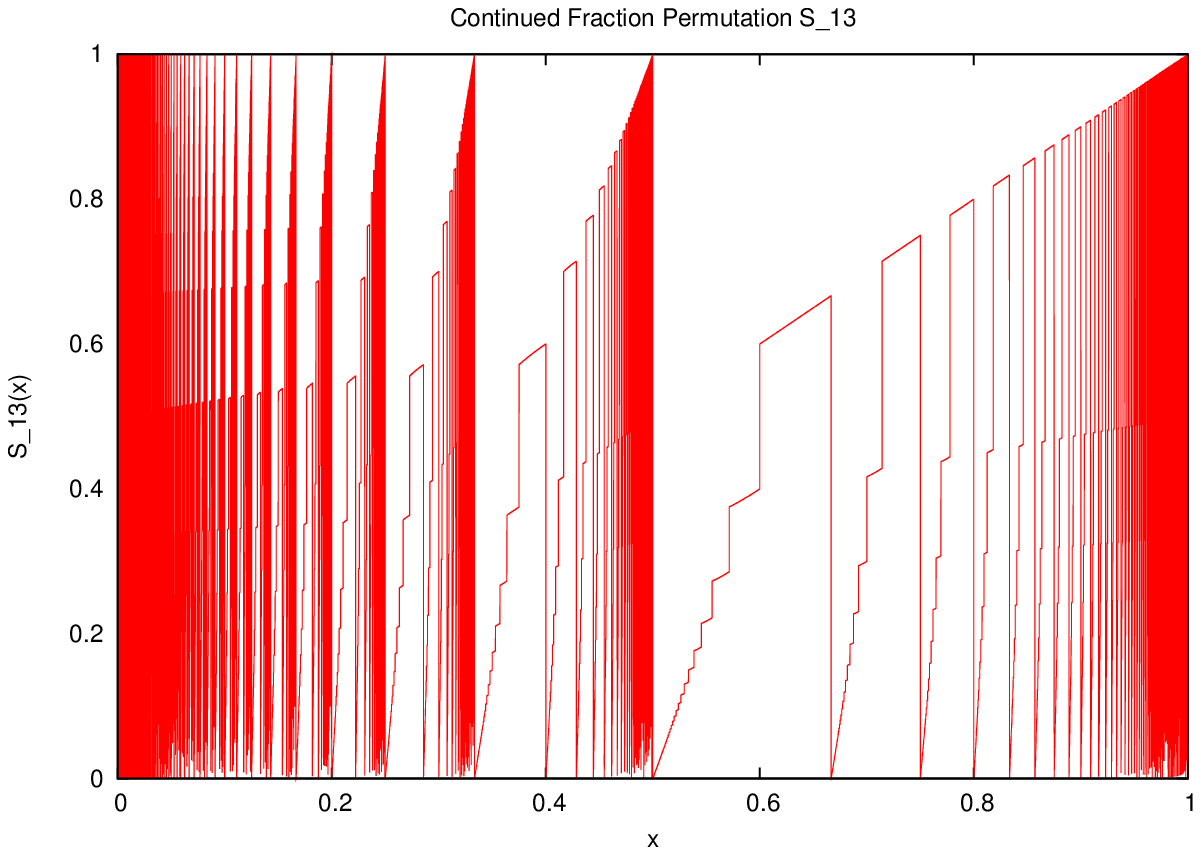}

\includegraphics[width=0.5\textwidth]{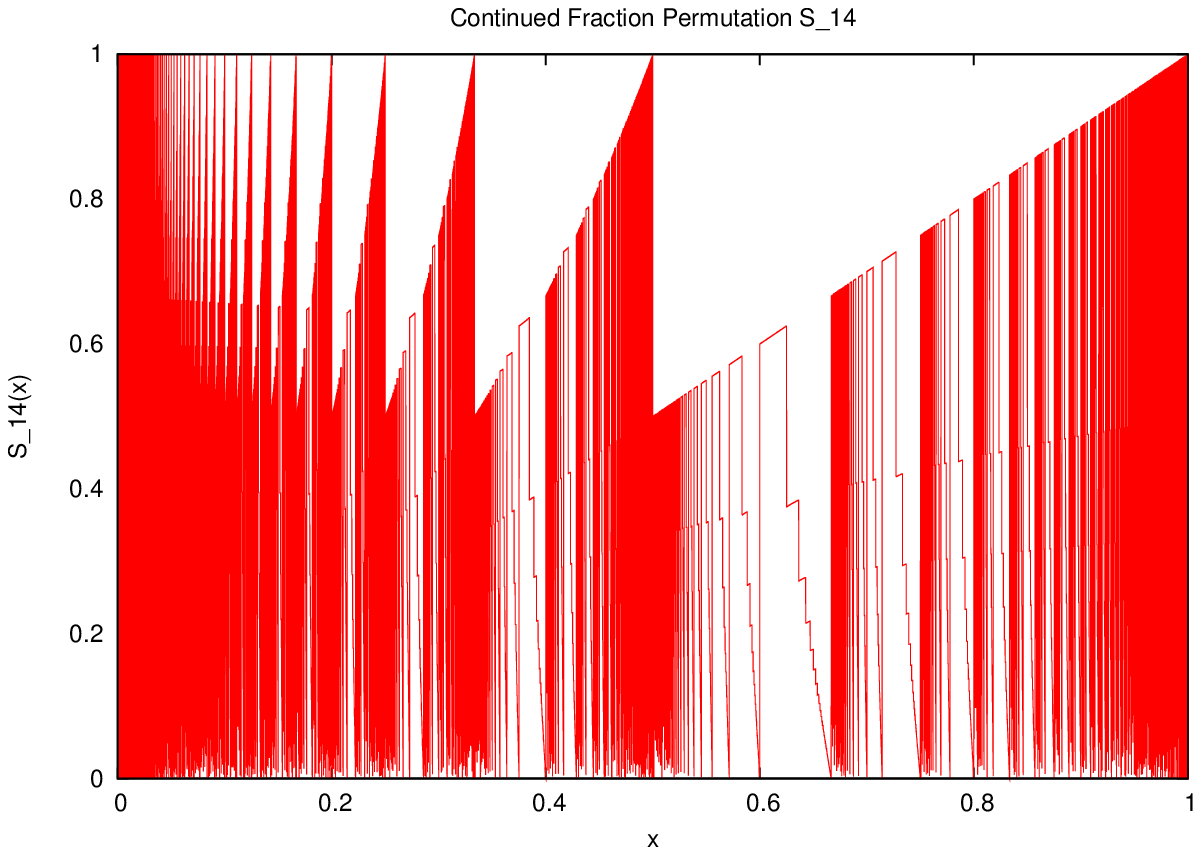}\includegraphics[width=0.5\textwidth]{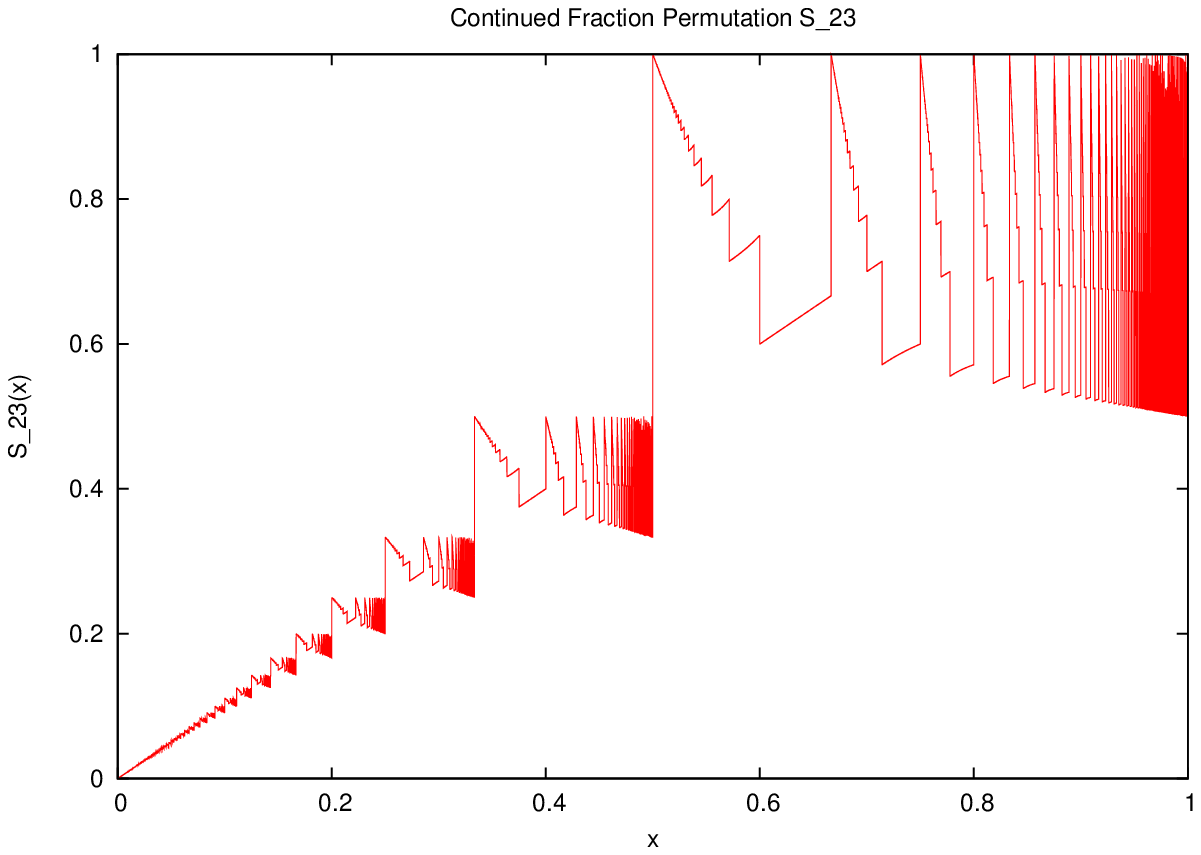}

The above figures show a set of graphs for various continued fraction
permutations $S_{p,q}(x)$ for the swap operator $S_{p,q}:[a_{1},\cdots,a_{p},\cdots,a_{q},\cdots]\mapsto[a_{1},\cdots,a_{q},\cdots,a_{p},\cdots]$.
A countable number of discontinuities are easily seen. It should also
be clear that the density of these discontinuities make it quite difficult
to work with these functions using simple-minded numerical approaches.
These figures posses a fractal self-similarity; this self-similarity
is presumably key to obtaining workable analytic expressions. 

\lyxline{\normalsize}
\end{figure}

Before exploring the sums \ref{eq:zeta-pq-sum} in detail, it is worth
examining the corresponding sum for the Riemann zeta, so as to have
a baseline to compare to. Define the Riemann sum as \[
\zeta_{N}\left(s\right)=\frac{s}{s-1}-\frac{s}{N}\;\sum_{n=0}^{N-1}h\left(\frac{2n+1}{2N}\right)\;\left(\frac{2n+1}{2N}\right)^{s-1}\]

\begin{figure}
\caption{\label{cap:Riemann-Sums}Riemann Sums}

\includegraphics[width=0.19\textwidth]{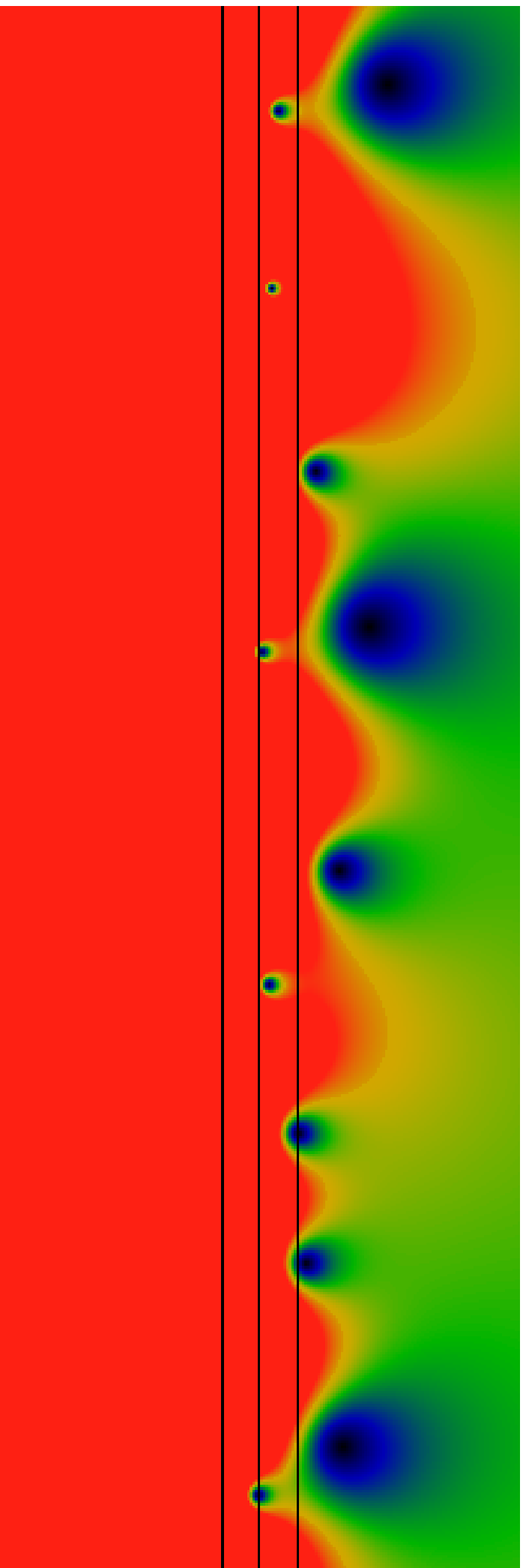} \includegraphics[width=0.19\textwidth]{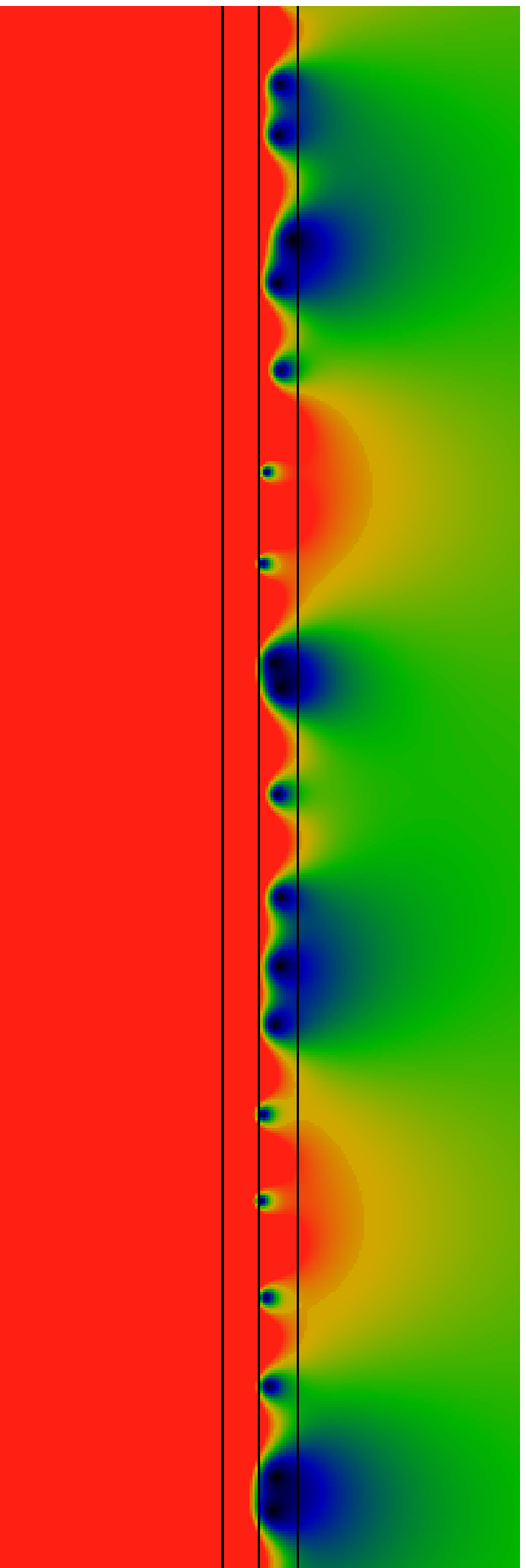}
\includegraphics[width=0.19\textwidth]{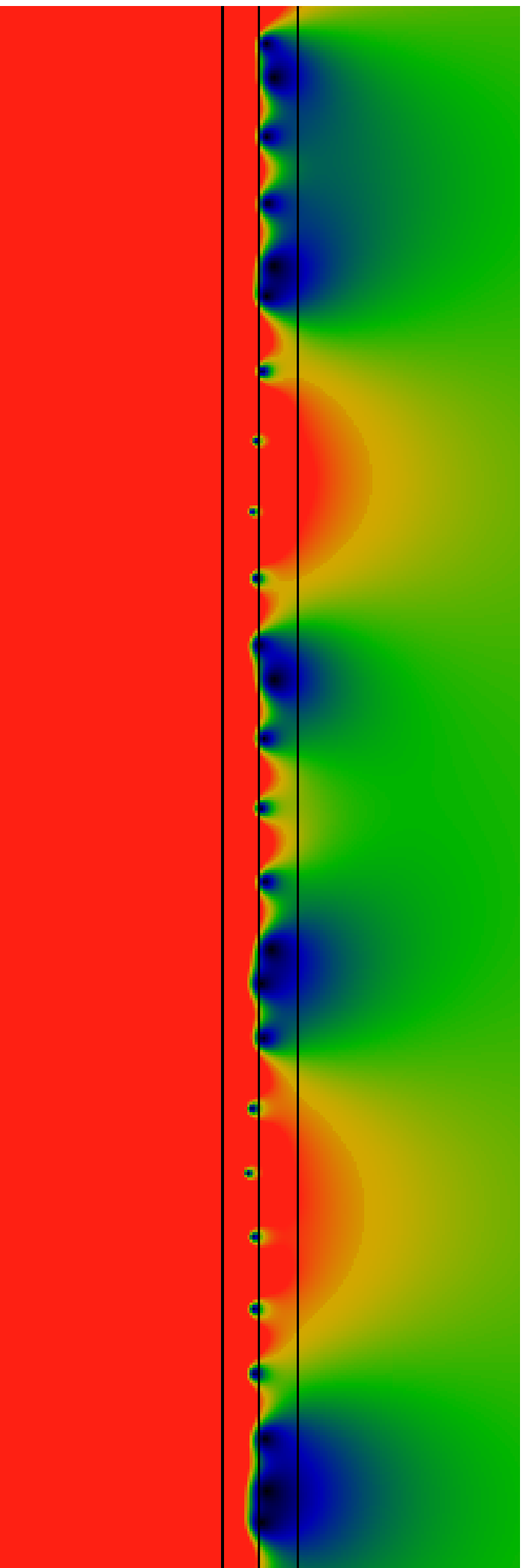} \includegraphics[width=0.19\textwidth]{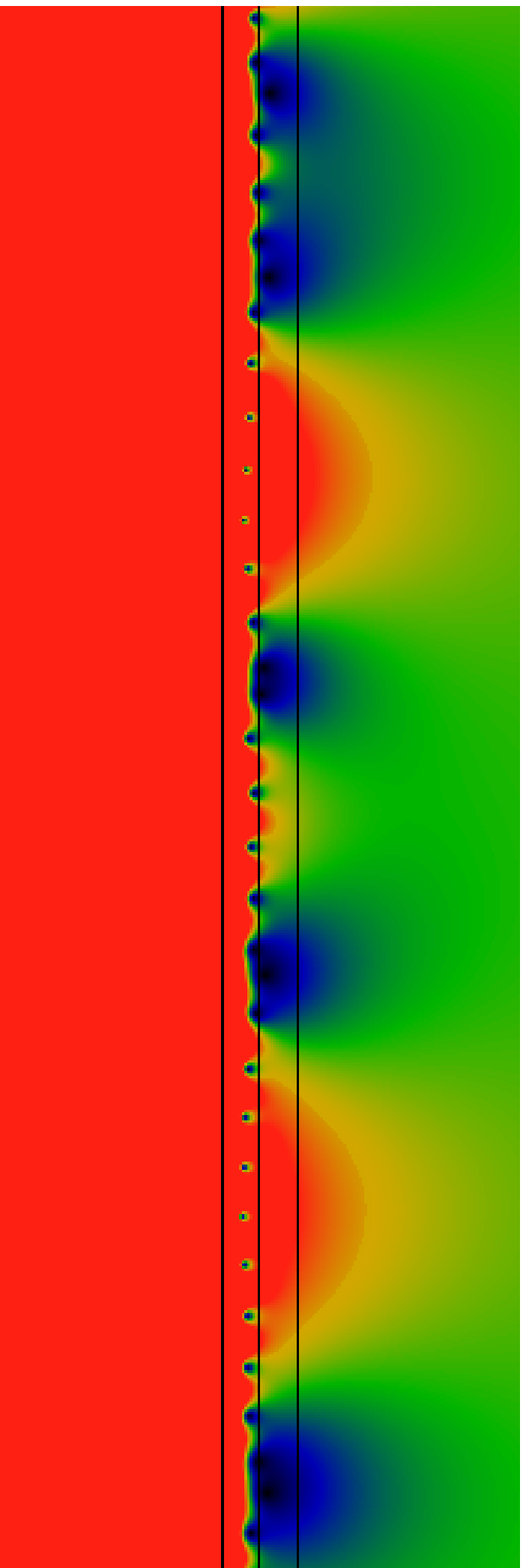}
\includegraphics[width=0.19\textwidth]{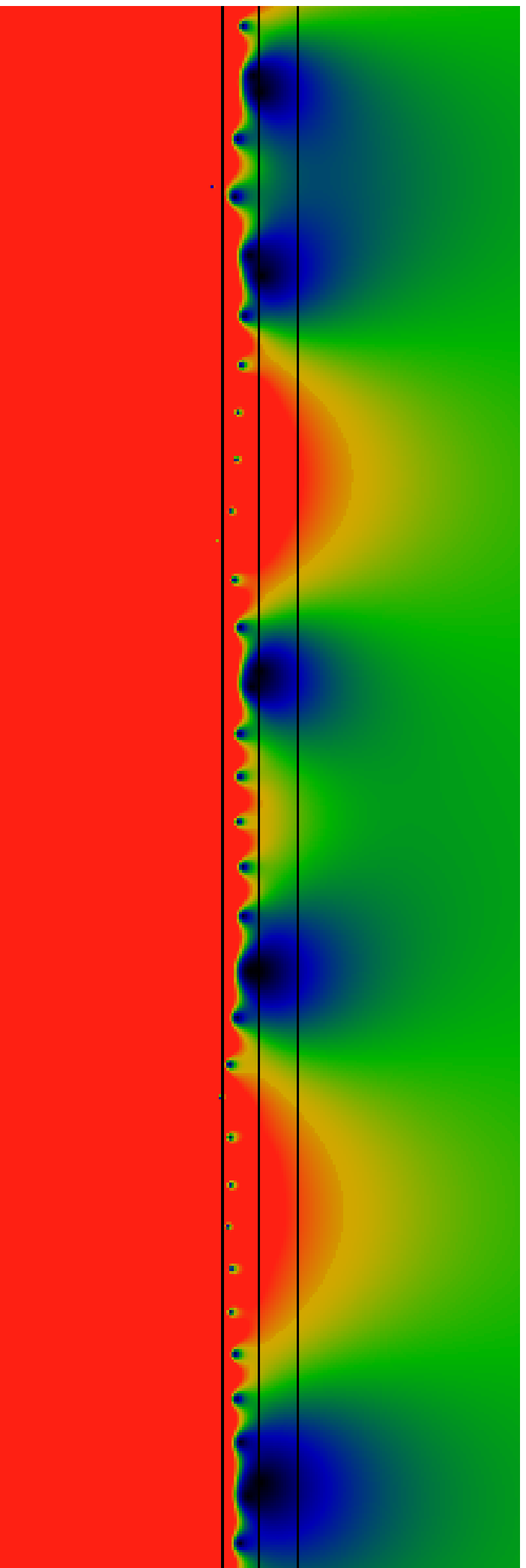}

The above strip charts show $\left|\zeta_{N}(s)\right|$ for N=15,
151, 1051, 11051 and 151051 respectively. The three vertical black
lines are located at $\Re s$=0, 1/2 and 1. The vertical range of
these strips are all identical, ranging from $\Im s$=13 at the bottom,
to $\Im s$=34 at the top. The coloration is such that red represents
values greater than two, green for values about equal to one, and
blue/black for values near/at zero. In this particular strip, we expect
the exact Riemann zeta to have five zeros, near 14.13, 21.02, 25.01,
30.42, and 32.93. Several remarkable properties are visible: First,
there are far more than five zeros visible, some lying outside the
critical strip. Next, as the value of $N$ increases, the zeros migrate
towards the $\Re s=1/2$ line; however, the number of zeros seems
to multiply logarithmically as well. Most of these zeroes seem to
have a tiny residue, which diminishes as $N$ increases; presumably,
they will completely disappear in the limit. What's left are five
zeros which are not fading away; these become the Riemann zeroes.

\lyxline{\normalsize}
\end{figure}

As a function of $N$, this sum converges very slowly to the Riemann
zeta function. As $N$ increases, the traditional zeros become visible,
although, or finite $N$, they do not sit exactly on the $\Re s=\frac{1}{2}$
critical line. The traditional zeros are also accompanied by a larger
number of {}``artifact zeros'': closely-spaced zeros, with a much
smaller residue, whose spacing (and residue) decreases as $N$ rises.
Visually, these are very easy to distinguish; the figure \ref{cap:Riemann-Sums}
show both the traditional and the artifact zeros clearly.

\begin{figure}
\caption{Riemann sum phase plots\label{fig:Riemann-sum-phase}}

\includegraphics[width=0.5\textwidth]{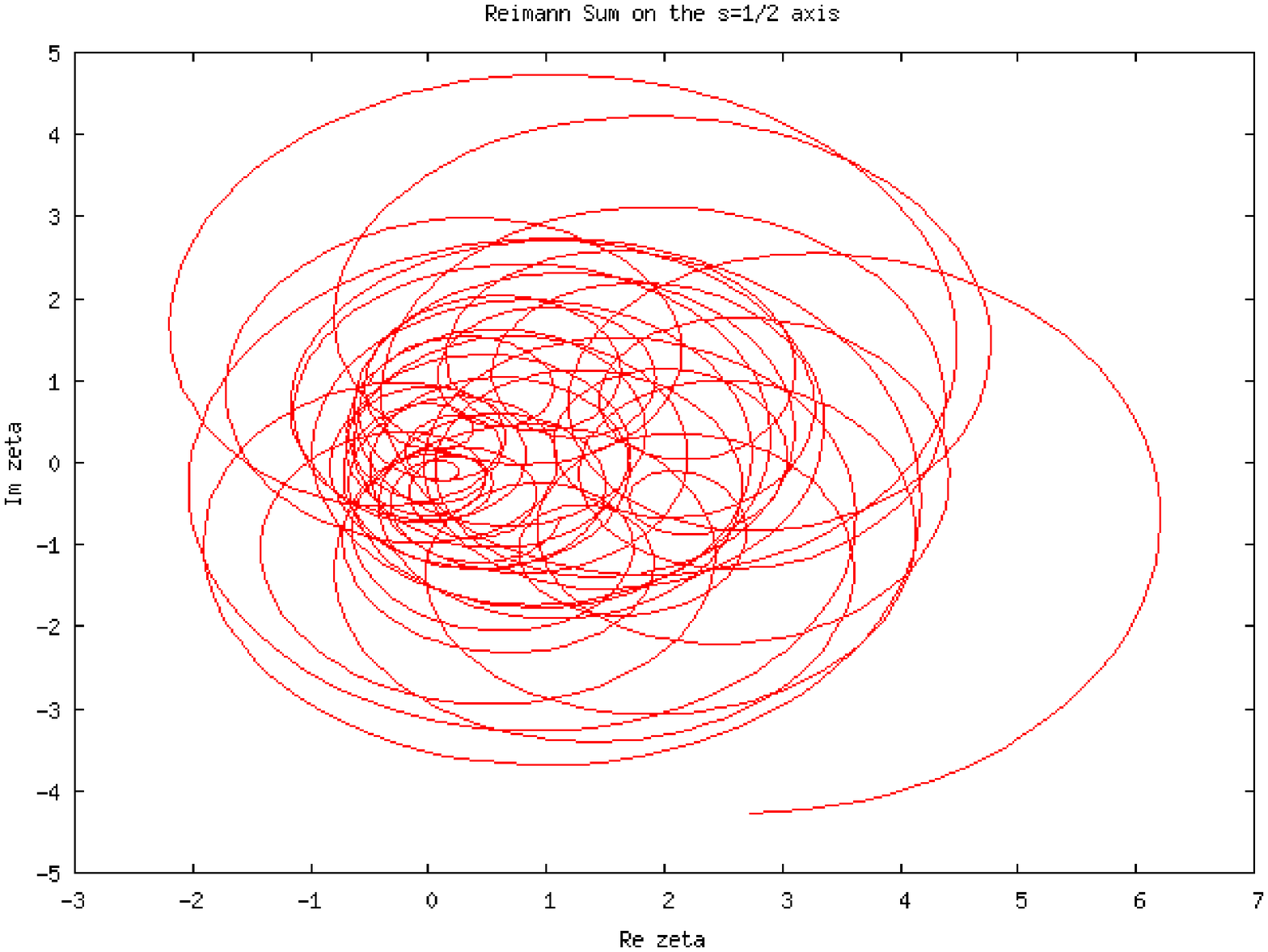}\includegraphics[width=0.5\textwidth]{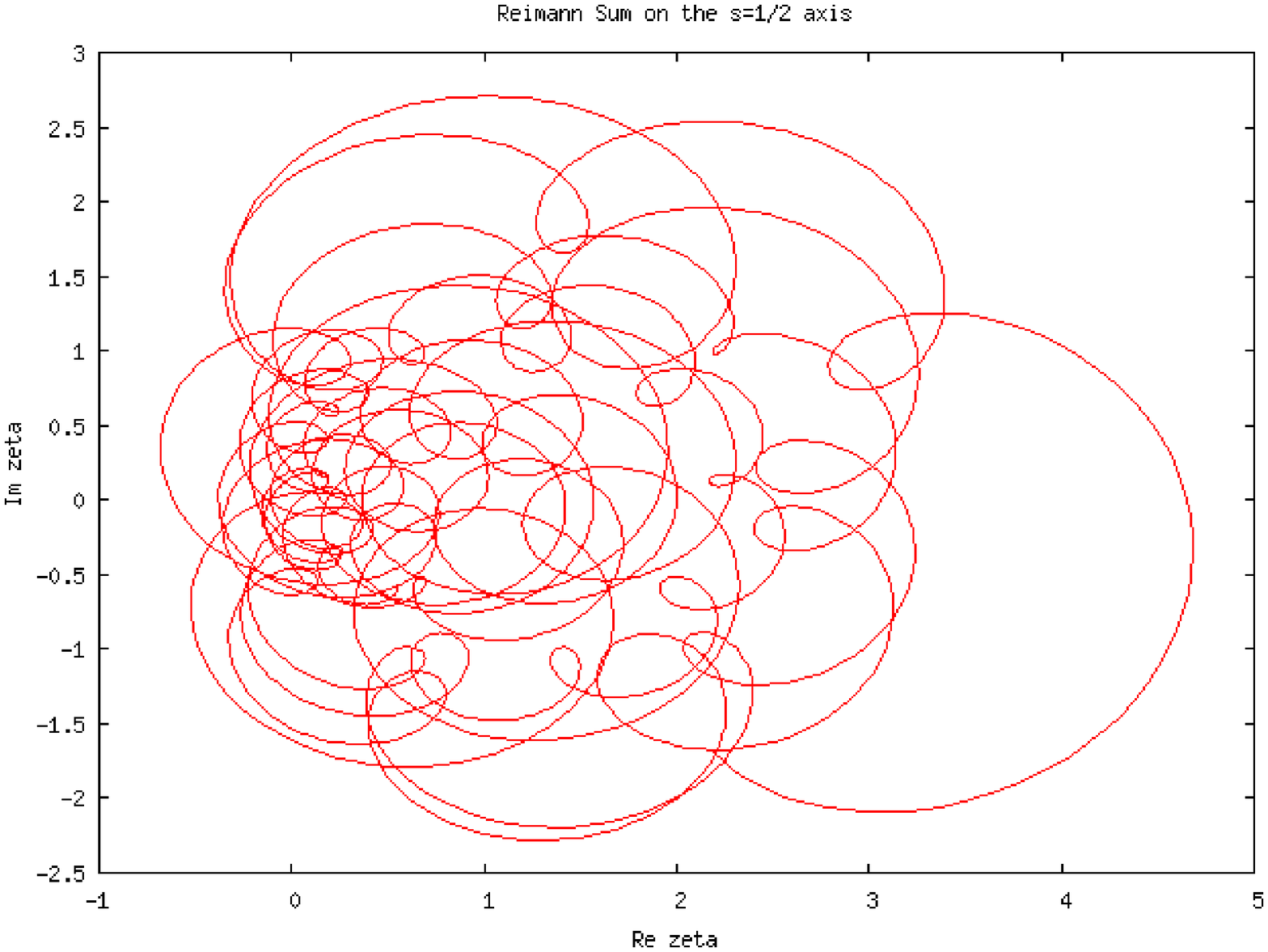}

\includegraphics[width=0.5\textwidth]{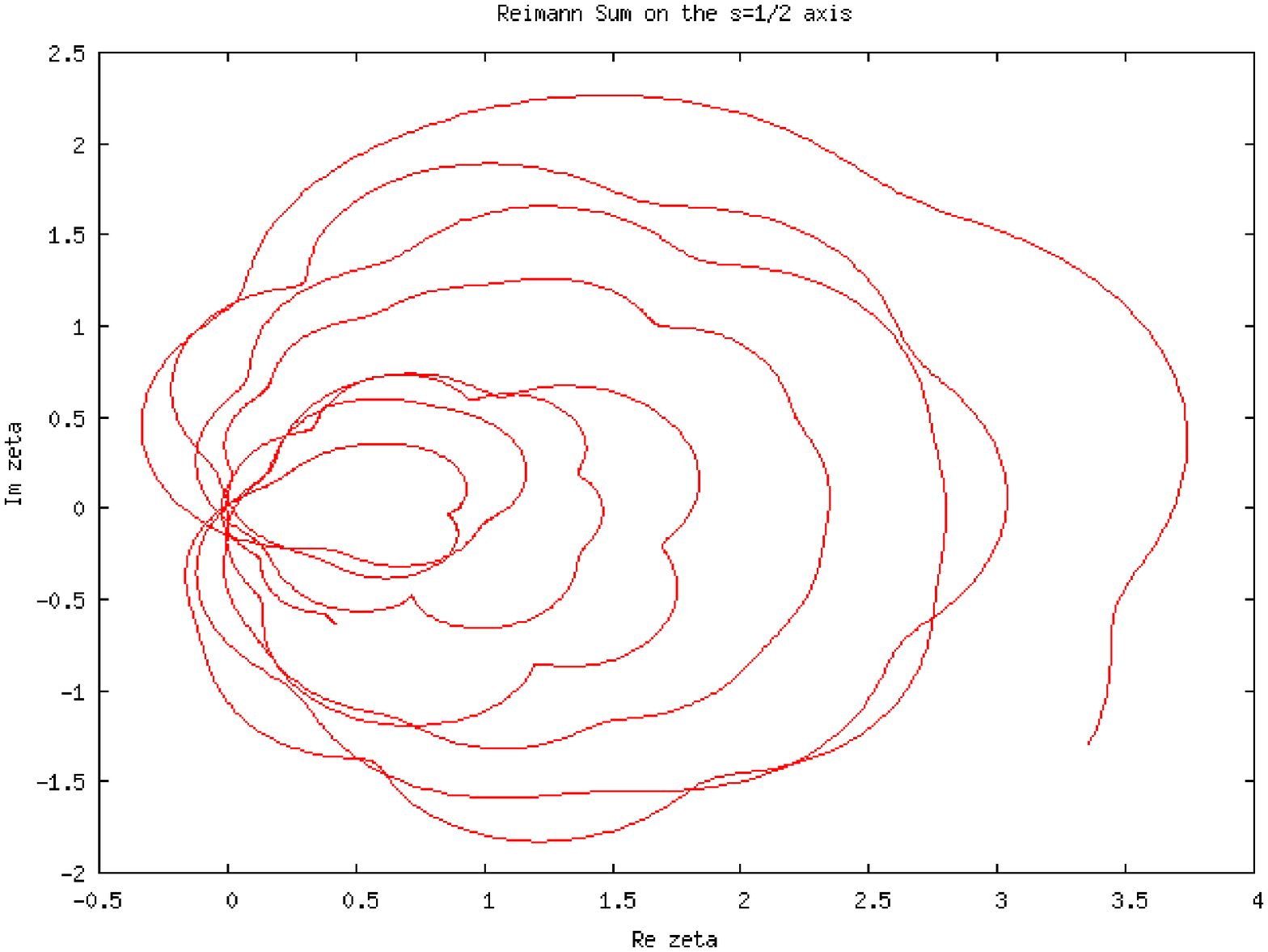}

The above three figures are parametric phase plots of $\zeta_{N}(s)$
for N=1051, 11051 and 121051, treating $s$ as a parameter, varying
from $\Im s$=13 to 34 along the $\Re s=1/2$ line. That is, the graphs
show $\Re\zeta_{N}$ on the x-axis, and $\Im\zeta_{N}$ along the
y-axis as $s$ varies. The large number of zeros in the critical strip
can be understood as being due to the looping seen in these graphs.
As $N$ increases, the loops shrink, passing through cusps; the last
figure beginning to resemble the true phase portrait of the Riemann
zeta.

\lyxline{\normalsize}
\end{figure}

The convergence properties of this sum, as well as the nature of the
{}``artifact zeros'', can be understood by exploring the phase of
$\zeta_{N}$ along the critical axis. This is shown in figure \ref{fig:Riemann-sum-phase}.
Truncation at finite $N$ introduces a quasi-sinusoidal perturbation
overlaid on the {}``true'' Riemann zeta function.

\begin{figure}
\caption{Zeros of $S_{1,2}$\label{fig:Zeros-of-S_12}}

\includegraphics[width=0.49\textwidth]{rzero-151051} \includegraphics[width=0.49\textwidth]{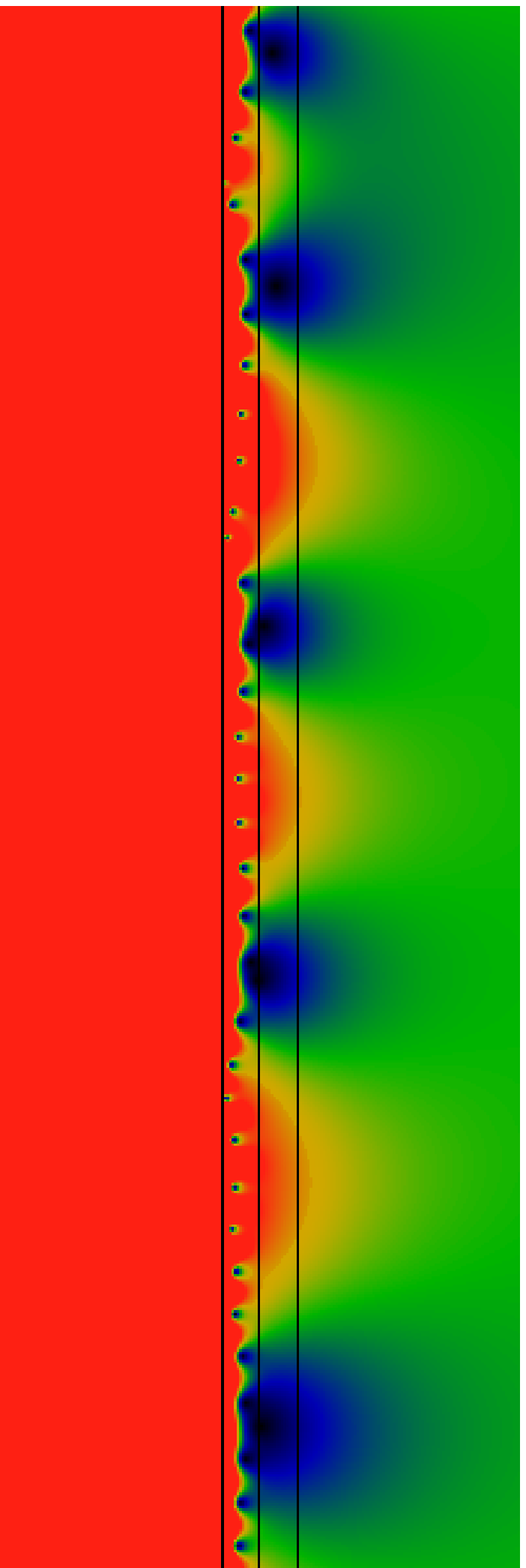}

This figure shows a side-by-side comparison of $\left|\zeta_{N}(s)\right|$
and $\left|\zeta_{N;1,2}(s)\right|$, with the later as defined in
equation\ref{eq:zeta-pq-sum}, and $N=151051$. The sizes, heights
and color scheme are as in the previous figures. Note the near-correspondence
of the dominant zeros; yet, the zeros for $\left|\zeta_{N;1,2}(s)\right|$
do appear to be at shifted locations.

\lyxline{\normalsize}
\end{figure}

\begin{figure}
\caption{Zeros of $S_{1,3}$ and $S_{2,3}$\label{fig:Zeros-of-S_13 and S_23}}

\includegraphics[width=0.49\textwidth]{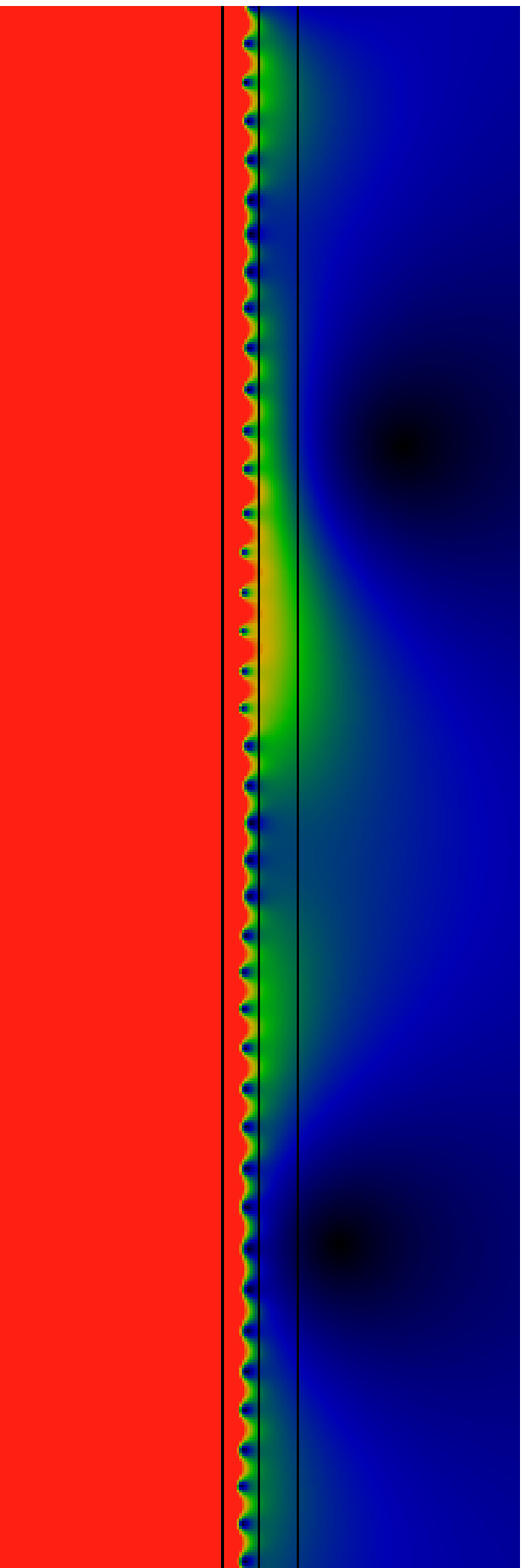} \includegraphics[width=0.49\textwidth]{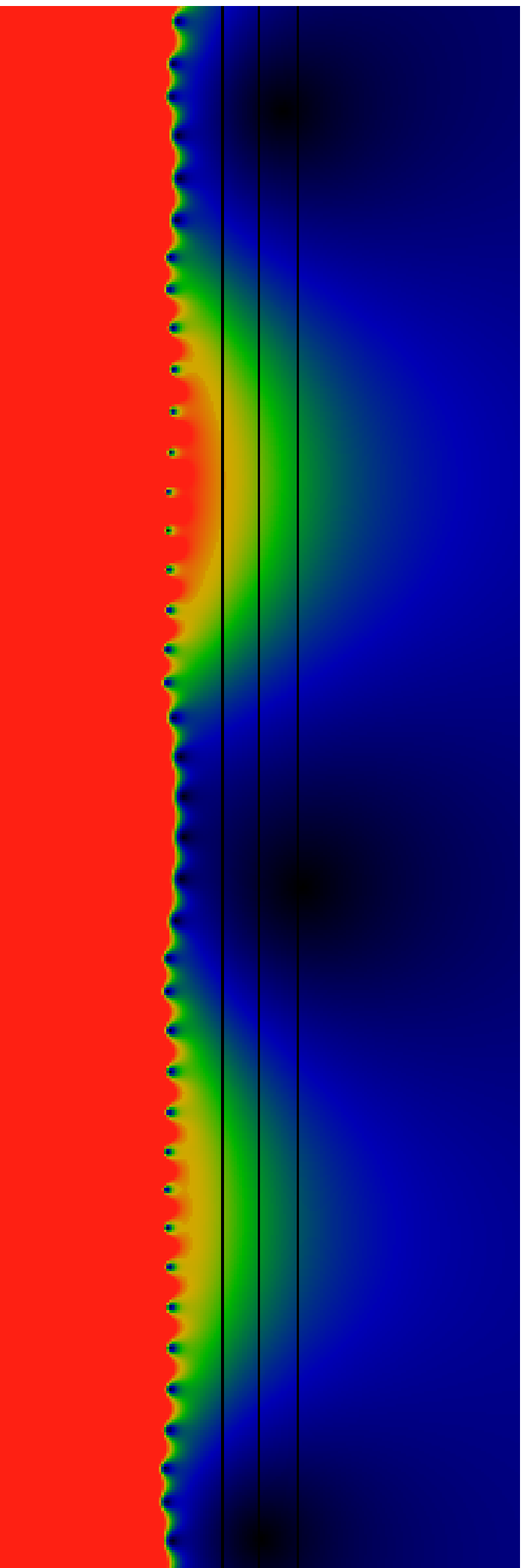}

This figure shows a side-by-side comparison of $\left|\zeta_{N;1,3}(s)\right|$
and $\left|\zeta_{N;2,3}(s)\right|$, with the later as defined in
equation\ref{eq:zeta-pq-sum}, and $N=151051$. The sizes, heights
and color scheme are as in the previous figures. The dominant zeros
no longer resemble those of $\left|\zeta_{N}(s)\right|$, and it is
no longer clear whether they will end up in the critical strip. However,
it is also the case that the structure of $S_{1,3}(x)$ and $S_{2,3}(x)$
is much finer and detailed, and so this value of $N$ does not yet
provide an adequate visual estimate of the $N\to\infty$ form of these
figures.

\lyxline{\normalsize}
\end{figure}

Armed with this overview, one may now explore the sum $\zeta_{N;1,2}$.
As figure \ref{fig:Zeros-of-S_12} shows, the simplest numerical evidence
seems to directly support the hypothesis. The evidence for the hypothesis,
in the case of $\zeta_{N;1,3}$ and $\zeta_{N;2,3}$ is less clear,
as seen in figure \ref{fig:Zeros-of-S_13 and S_23}. On the other
hand, the numeric sums, for any fixed $N$, are far less {}``accurate''
than that for $\zeta_{N;1,2}$. That is, the discontinuities for $S_{1,3}$
and $S_{2,3}$ are very finely spaced, and the summation does a considerably
poorer job of capturing these. 

\section{A Remarkable Shadow}

The numerical evaluation of sums requires a certain amount of cross-checking
so as to catch errors. A particularly curious sum to evaluate is\[
\eta_{N}(s)=\frac{s}{s-1}-\frac{s}{N}\;\sum_{n=0}^{N-1}\left(\frac{2n+1}{2N}\right)^{s}\]
which does no permutation at all. That is, one effectively has $\eta_{N}(s)=\zeta_{N;p,p}(s)$
for any integer $p$. The limiting integral is trivial to calculate,
since $S_{p,p}(x)=x$, and so\[
\lim_{N\to\infty}\eta_{N}(s)=\frac{2s}{s^{2}-1}\]
Thus, this sum shouldn't be interesting, except as a calibration of
sorts. 

\begin{figure}
\caption{\label{cap:A-Shadow-of}A Shadow of Zeta}

\includegraphics[width=0.49\textwidth]{rzero-11051} \includegraphics[width=0.49\textwidth]{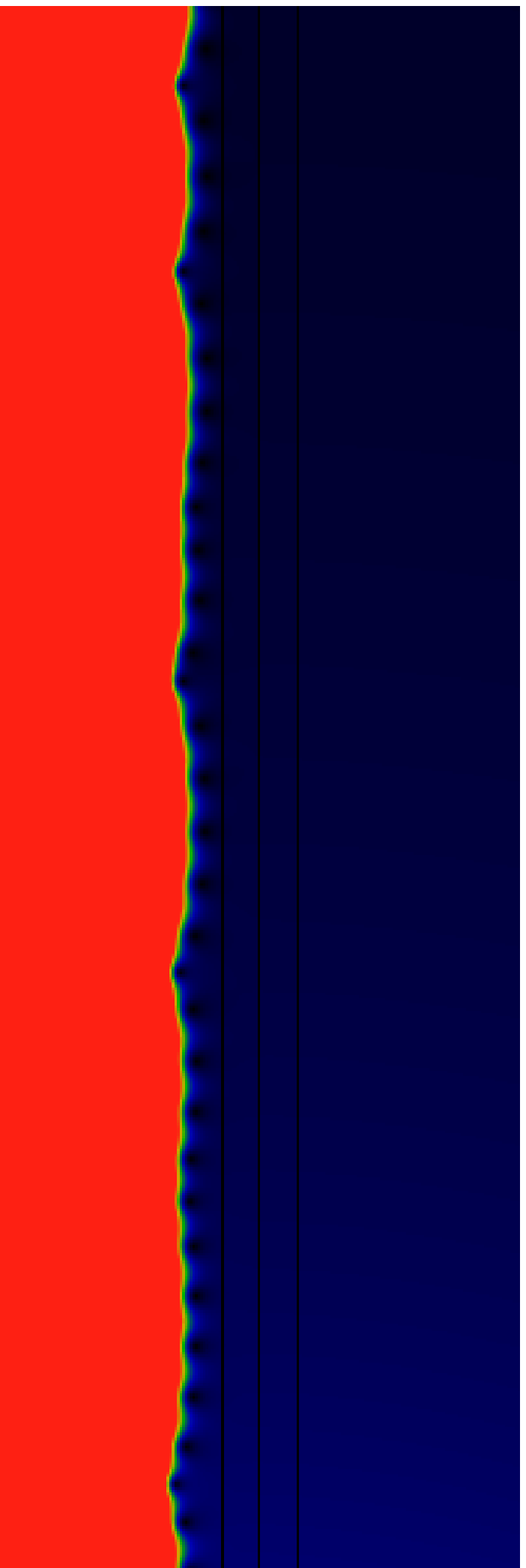}

This figure compares $\left|\zeta_{N}(s)\right|$ and $\left|\eta_{N}(s)\right|$
for $N=11051$. Note the {}``horns'' or {}``cusps'' in $\left|\eta_{N}(s)\right|$.
While these are falling well outside the critical strip, it is particularly
remarkable that they line up vertically with the locations of the
non-trivial Riemann zeros.

\lyxline{\normalsize}
\end{figure}

But that is false. A review of the figure \ref{cap:A-Shadow-of} shows
that this simple stair-step sum retains traces of the locations of
the non-trivial Riemann zeros!

The presence of this shadow is not hard to explain. Define\[
Z_{N}(s)=\sum_{n=1}^{N}\frac{1}{n^{s}}\]
Then it is clear that, for large $N$, $Z_{N}(s)$will approximate
the Riemann zeta function:\[
\lim_{N\to\infty}Z_{N}(s)=\zeta(s)\]
at least for the case of $\Re s>1$, where the sum is convergent.
With some simple manipulations, one may rewrite the numeric sum as
\[
\sum_{n=0}^{N-1}\left(\frac{2n+1}{2N}\right)^{s}=\frac{1}{(2N)^{s}}\left[Z_{2N-1}(-s)-2^{s}Z_{N-1}(-s)\right]\]
The fact that the seemingly {}``trivial'' $\eta_{N}(s)$ sum contains
this approximation to the zeta would seem to explain the numeric features
shown in figure \ref{cap:A-Shadow-of}. Furthermore, the presence
of the $-s$ in the arguments appears to explain why the horns are
located not far from the $\Re s=-1/2$ line. 

On the other hand, there is another important lesson to be drawn here:
this {}``trivial'' sum is insufficient to explain the features seen
in the $\zeta_{N;p,q}(s)$ sums. That is, although the numerical evaluation
itself is prone to introducing the Riemann zeta as an artifact, the
strength of this artifact does not immediately appear to be dominant
enough to account for the entire behavior of the $\zeta_{N;p,q}(s)$
sums.

\section{Sums and Integrals\label{sec:Sums-and-Integrals}}

This section explores analytic approaches to evaluating the integral
\ref{eq:zeta-permute}. Because $S_{p,q}\left(x\right)$ is a piecewise
combination of Mobius functions \[
\frac{ax+b}{cx+d}\]
it is enough to evaluate the indefinite integral 

\[
\int\frac{ax+b}{cx+d}x^{s-1}dx\]
This integral does not have any simple, closed-form solution, but
can be given in terms of an infinite sum. This may be obtained via
the Newton series \[
\left(1+x\right)^{s}=\sum_{k=0}^{\infty}\dbinom{s}{k}x^{k}\]
where \[
\binom{s}{k}=\frac{\Gamma\left(s+1\right)}{\Gamma\left(k+1\right)\Gamma\left(s-k+1\right)}\]
is the binomial coefficient. The series is absolutely convergent on
the unit interval $0\le x\le1$ that we wish to explore. The indefinite
integral may then be evaluated as a sum:\begin{align*}
F\left(s,\alpha;x\right)= & \int\frac{x^{s}}{x+\alpha}dx\\
= & \int\frac{\left(y-\alpha\right)^{s}}{y}dy\\
= & \left(-\alpha\right)^{s}\sum_{k=0}^{\infty}\binom{s}{k}\left(\frac{-1}{\alpha}\right)^{k}\int y^{k-1}dy\\
= & \left(-\alpha\right)^{s}\left[\ln\left(x+\alpha\right)+\sum_{k=1}^{\infty}\frac{\left(-1\right)^{k}}{k}\binom{s}{k}\left(1+\frac{x}{\alpha}\right)^{k}\right]\end{align*}
This allows the Mobius integral to be written as\[
\int\frac{ax+b}{cx+d}x^{s-1}dx=\frac{a}{c}F\left(s,\frac{d}{c};x\right)+\frac{b}{c}F\left(s-1,\frac{d}{c};x\right)\]
Using that fact that, in general, for continued fractions, the $a,b,c$
and $d$ will be integers, with $ad-bc=\pm1$, one may re-write the
above as\begin{equation}
\int\frac{ax+b}{cx+d}x^{s-1}dx=\frac{\left(-d\right)^{s-1}}{c^{s+1}}\left[\mp\ln\left(x+\frac{d}{c}\right)+\sum_{k=1}^{\infty}\frac{\left(-1\right)^{k}}{k}\binom{s-1}{k}\left(1+\frac{cx}{d}\right)^{k}\left(\frac{\pm s+kbc}{k-s}\right)\right]\label{eq:mobius sums}\end{equation}

For $S_{1,2}\left(x\right)$, this can be made more explicit. So,
we have that \[
S_{1,2}\left(x\right)=\frac{1}{a_{2}+\frac{1}{a_{1}+r}}\]
when $x$ is written as \[
x=\frac{1}{a_{1}+\frac{1}{a_{2}+r}}\]
with $a_{1}=\left\lfloor 1/x\right\rfloor $ and $a_{2}=\left\lfloor x/\left(1-a_{1}x\right)\right\rfloor $
integers. Solving for $r$ and replacing, one gets\[
S_{1,2}\left(x\right)=\frac{ax+b}{cx+d}\]
with \begin{align*}
a= & 1+a_{1}\left(a_{2}-a_{1}\right)\\
b= & a_{1}-a_{2}\\
c= & \left(a_{2}-a_{1}\right)\left(1+a_{1}a_{2}\right)\\
d= & 1+a_{2}\left(a_{1}-a_{2}\right)\end{align*}
Note that in this case, $ad-bc=1$ for all possible values of $a_{1}$
and $a_{2}$. This form is valid over the interval \[
\left(\underline{x},\overline{x}\right)=\left(\frac{a_{2}}{1+a_{1}a_{2}},\frac{1+a_{2}}{1+a_{1}+a_{1}a_{2}}\right)\]
and so one may combine these elements to write \begin{equation}
\int_{0}^{1}S_{1,2}\left(x\right)x^{s-1}dx=\sum_{a_{1}=1}^{\infty}\sum_{a_{2}=1}^{\infty}\int_{\underline{x}}^{\overline{x}}\frac{ax+b}{cx+d}x^{s-1}dx\label{eq:analytic sum}\end{equation}

At this point, the ability to perform any further simple manipulations
stops. The summations are only conditionally convergent; one cannot
interchange the sum over $k$ with the sums over $a_{1}$, $a_{2}$.
If one were to do so, one would quickly discover that, for fixed $k$,
the sums over $a_{1}$, $a_{2}$ are logarithmically divergent. While
this is not immediately obvious just be gazing at the sums, it is
easily confirmed by numerically evaluating them. 

One more avenue suggests itself, but then fails: one might consider
applying the binomial expansion to the powers $\left(1+\frac{x}{\alpha}\right)^{k}$
that appear in the sum for $F\left(s,\alpha;x\right)$, and then exchanging
the order of summations there. But this is not possible; the sums
fail outright. Thus, it appears that there are no further (obvious)
analytic techniques that one can apply to reduce the expression \ref{eq:analytic sum}.
One is then left to contemplate the numeric evaluation of \ref{eq:analytic sum}.
Naively, one might expect that doing so would offer considerable advantages
over the numeric sum \ref{eq:zeta-pq-sum}, but this is not at all
the case; convergence is even slower, and fails more profoundly. 

Thus, the remainder of this text uses only the sum \ref{eq:zeta-pq-sum}
for numeric exploration. The position of the first non-trivial zero
of $\zeta_{1,2}\left(s\right)$ is examined in the next section.

\section{Locating the first zero}

This section explores the numeric evidence for the location of the
first zero of $\zeta_{1,2}\left(s\right)$. One may quickly discover
that it must be near $s=\frac{1}{2}+i\,14.92$, but improving upon
this value is remarkably difficult, given only the tools \ref{eq:zeta-pq-sum}
and \ref{eq:analytic sum}.

\begin{figure}
\caption{Numerical zero-finding\label{fig:Numerical-zero-finding}}

\includegraphics[width=0.5\textwidth]{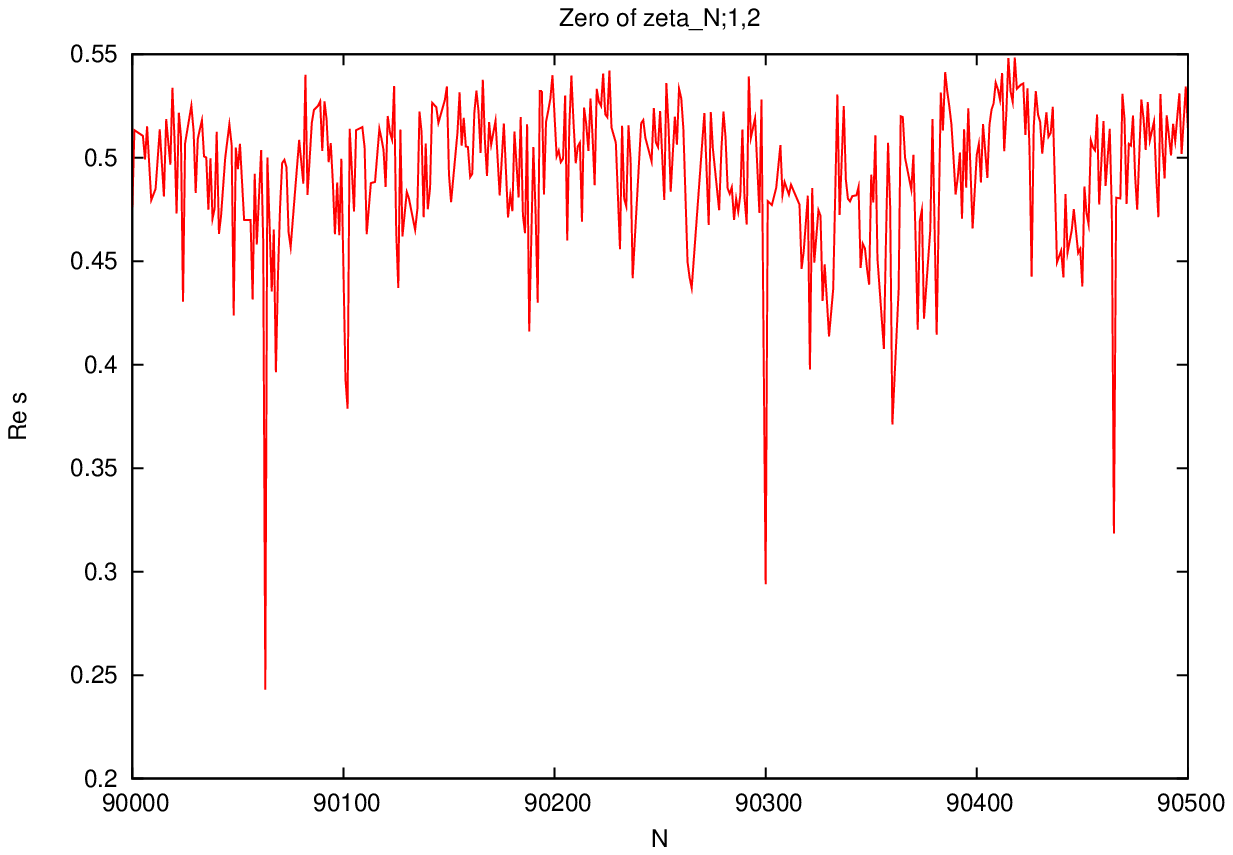}\includegraphics[width=0.5\textwidth]{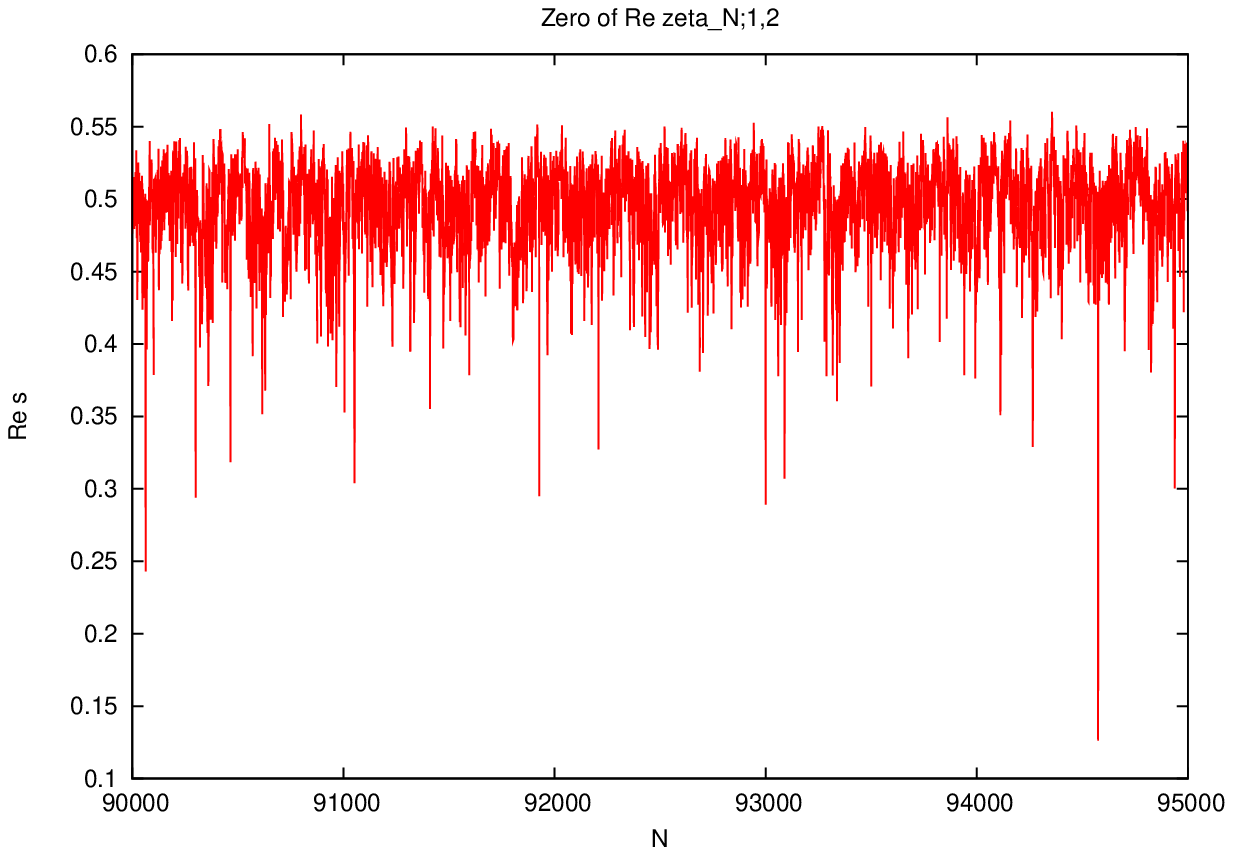}

\includegraphics[width=0.5\textwidth]{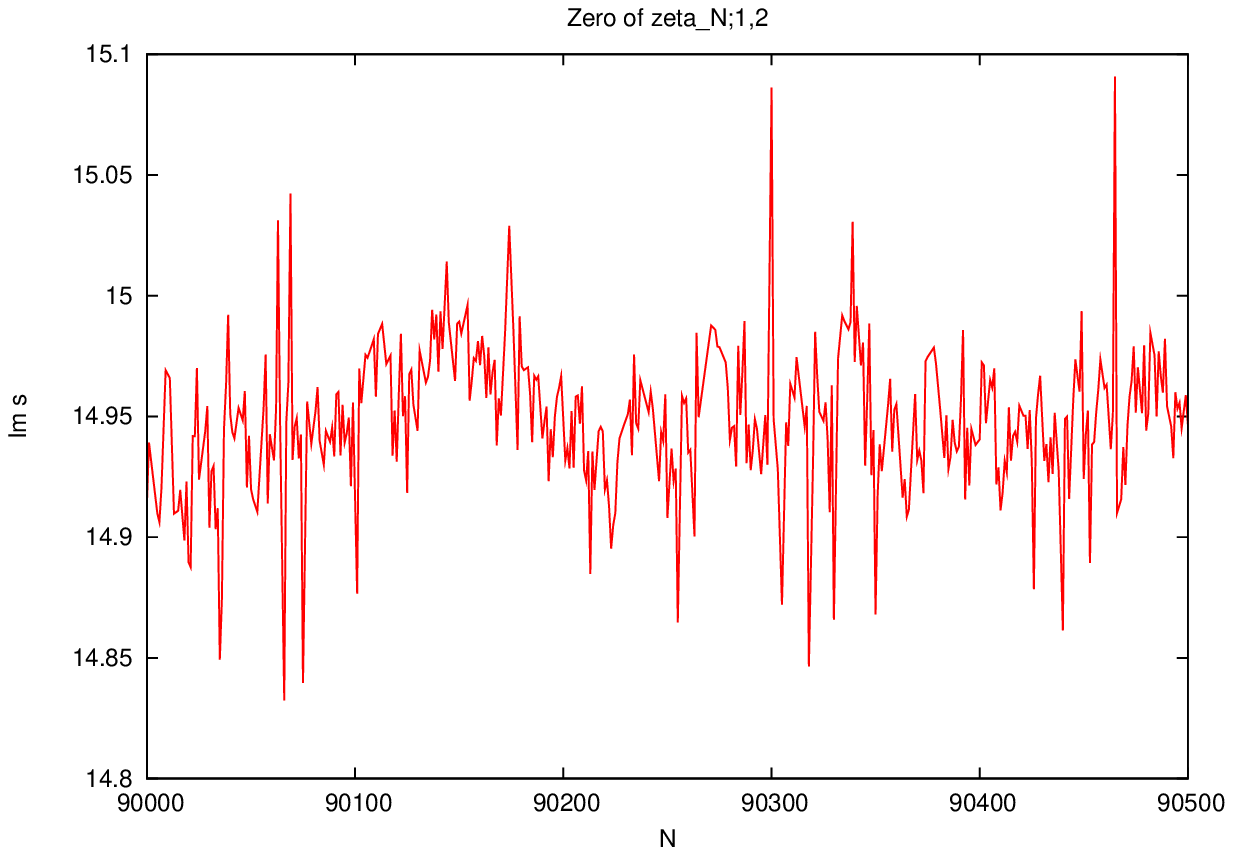}\includegraphics[width=0.5\textwidth]{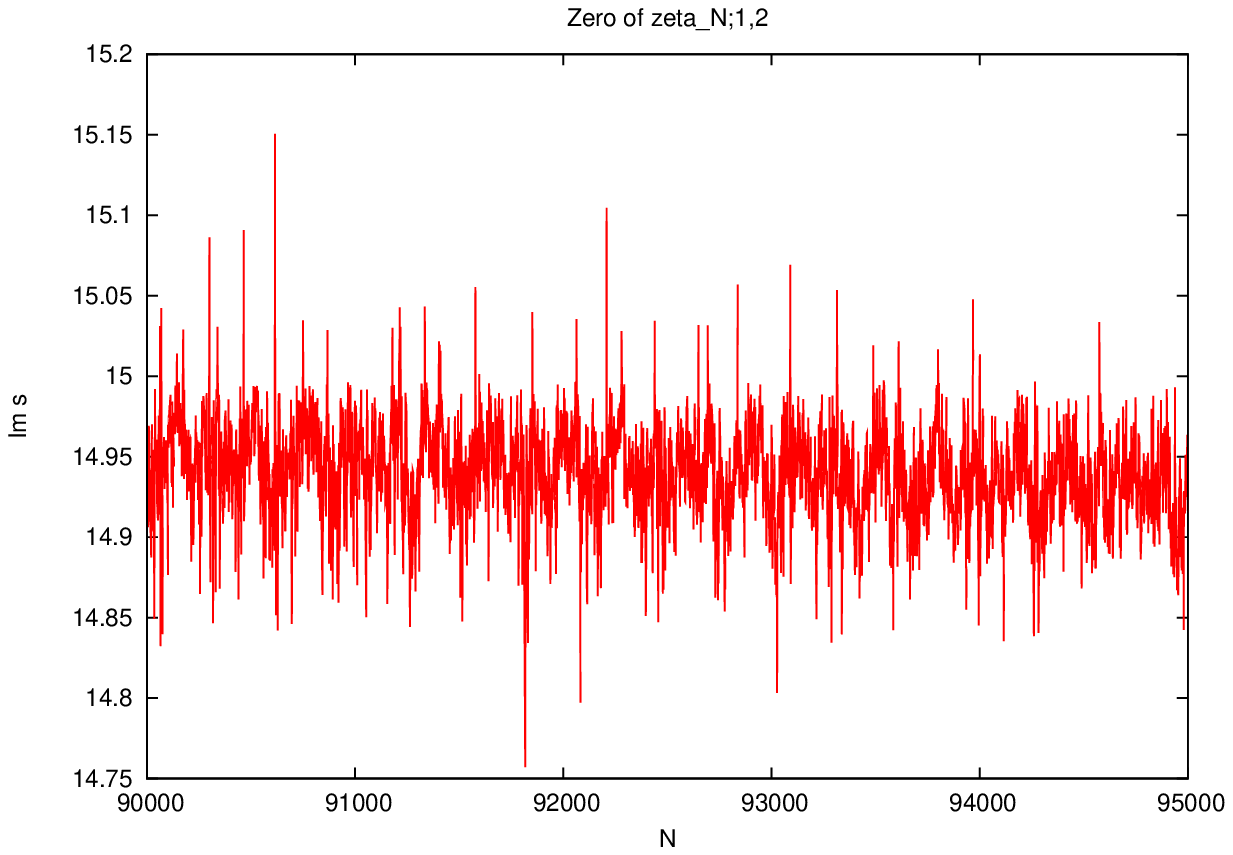}

The above graphs show the location of the first non-trivial zero of
$\zeta_{N;1,2}\left(s\right)$, as a function of $N$, in the range
of $90000\le N\le95000$. As may be seen, the location of the zero
depends very strongly on $N$, and is an extremely noisy function
of $N$. There is the vaguest hint of an oscillatory behavior in $N$,
possibly with a period of $\Delta N=300$, but this is very strongly
obscured by the noise. The power spectrum of this same data-set is
examined in the next figure. 

\lyxline{\normalsize}

\end{figure}

The summations \ref{eq:zeta-pq-sum} are straightforward to evaluate
numerically. Estimates for the location of the zero may be obtained
using more-or-less standard zero-finding techniques, such as Powell's
algorithm\cite{Pre88} (some care must be taken, as the neighborhood
of the zero appears to be a bit wobbly). The resulting estimates for
the zero depend very strongly on $N$, and are very noisy, as the
figure \ref{fig:Numerical-zero-finding} shows. The noise appears
to be classic fractal $1/f$ noise, as can be seen from the graph
of the power spectrum.

\begin{figure}
\caption{Power spectrum}

\includegraphics[width=0.5\textwidth]{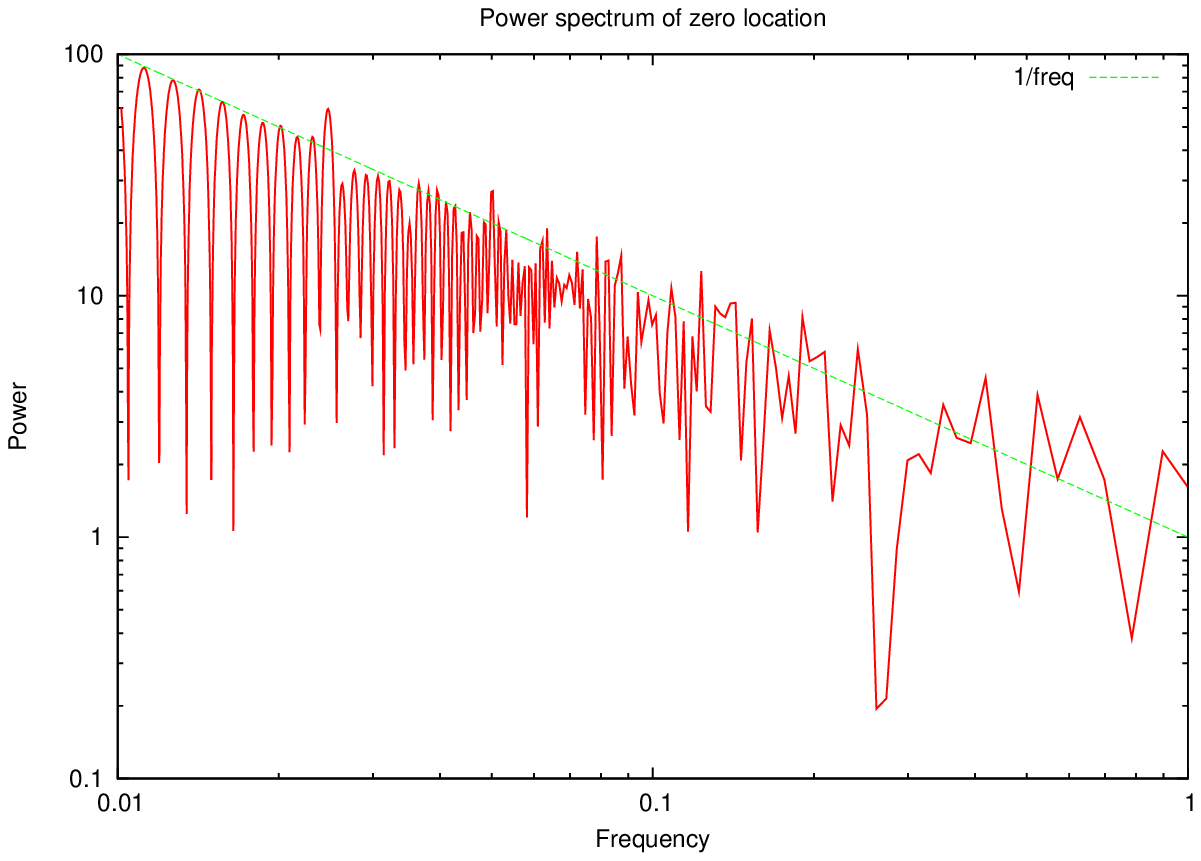}\includegraphics[width=0.5\textwidth]{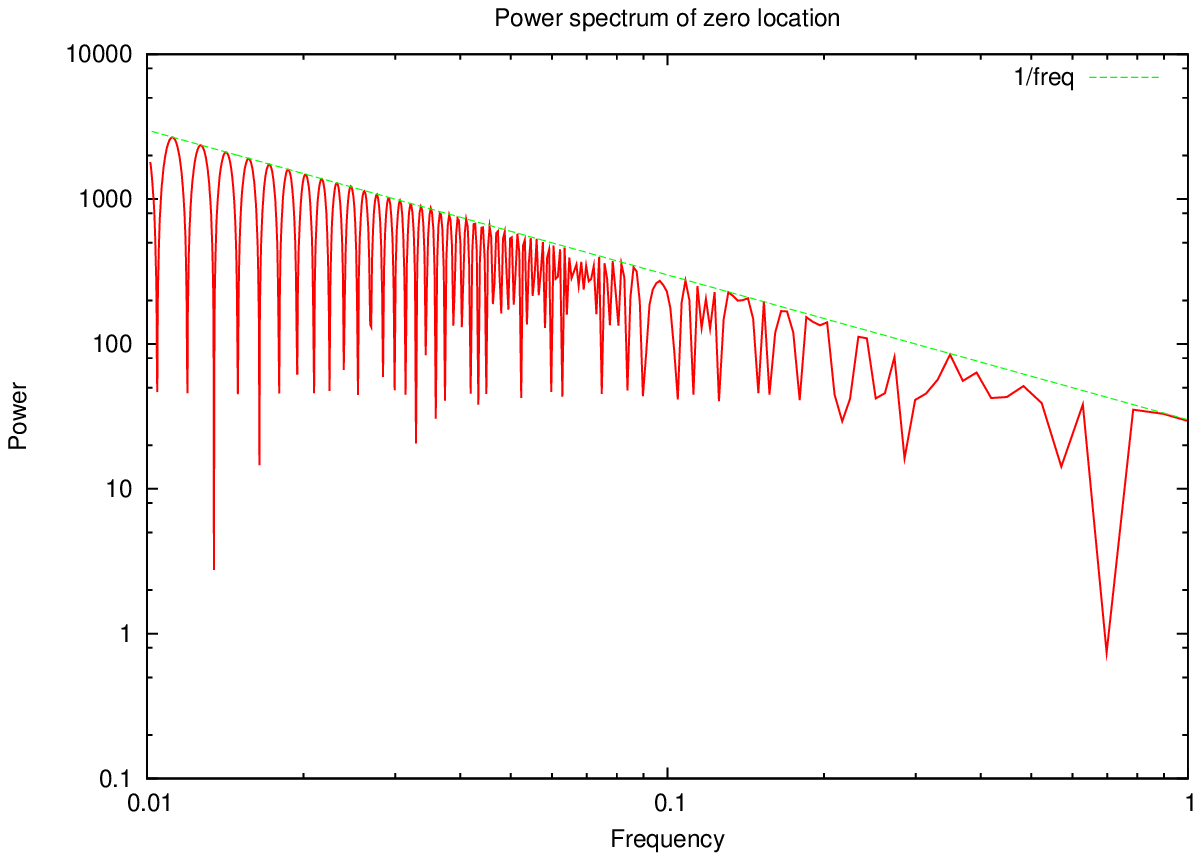}

These two charts show the power spectrum of the zero locations graphed
in figure \ref{fig:Numerical-zero-finding}. The power spectrum is
defined as usual, it is the square of the Fourier transform of the
data. The straight line shows the $1/f$ spectrum; the noise in the
location of the zero clearly obeys the classic fractal $1/f$ power
density. The regular pattern of arches appearing at low frequencies
is a numerical artifact resulting from the limited size of the data
sample. The left chart shows the power along the real axis (resulting
in the uncertainty of the $\Re s=\frac{1}{2}$ location), while the
right chart shows the power in the imaginary direction. The higher
power in the right chart indicates that the location $\Im s\approx14.92$
is more certain than the location $\Re s\approx\frac{1}{2}$. Curiously,
the {}``more important'' location is harder to pin down! 

\lyxline{\normalsize}

\end{figure}

\begin{figure}

\caption{Increased precision}
 \includegraphics[width=0.5\textwidth]{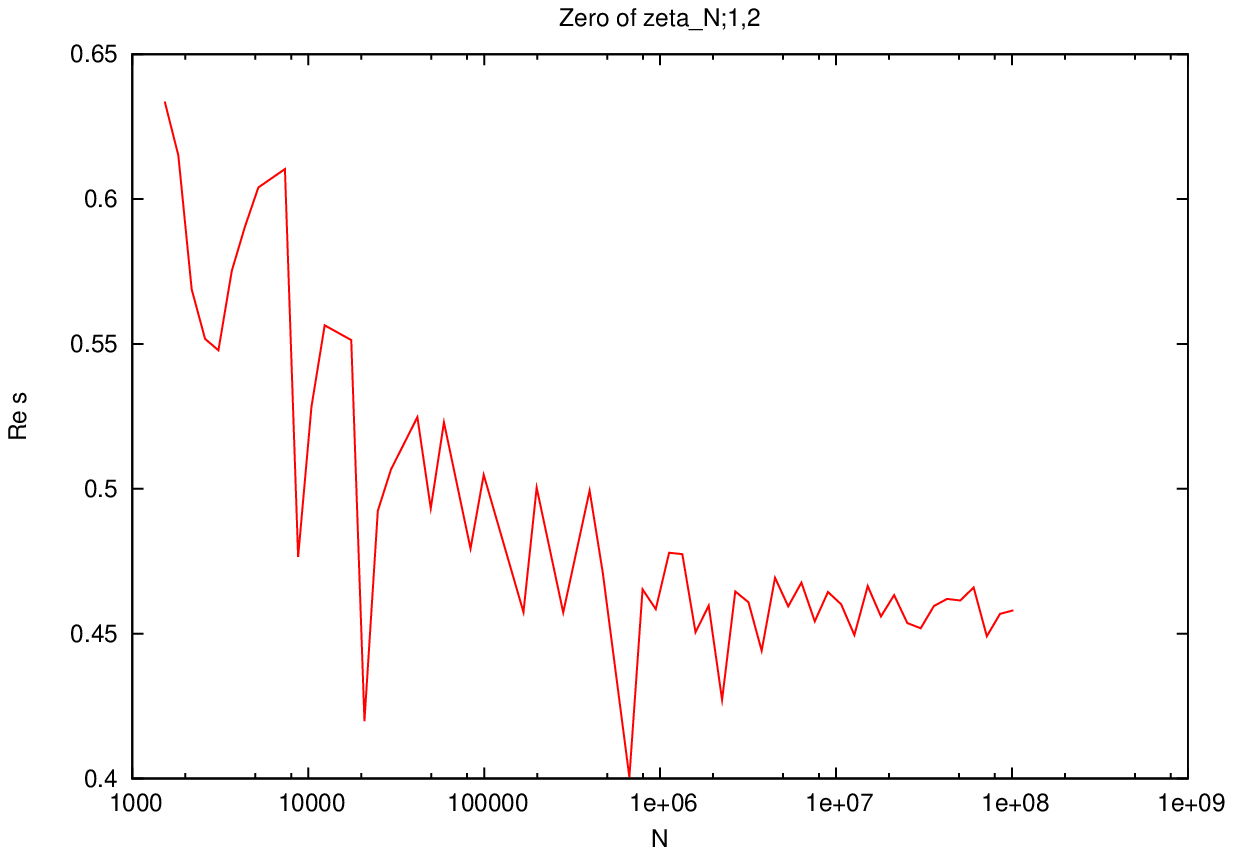}\includegraphics[width=0.5\textwidth]{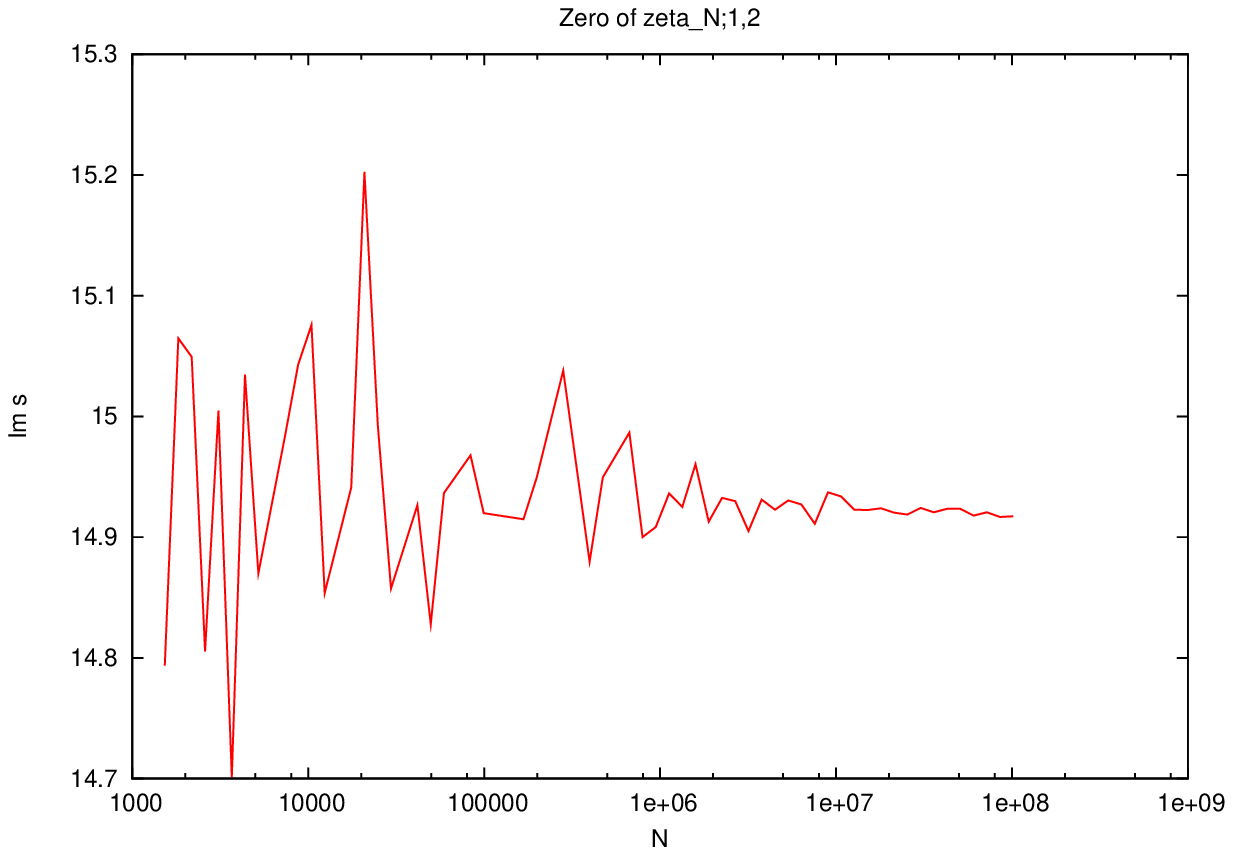}

The above graphs show the location of the first non-trivial zero of
$\zeta_{N;1,2}\left(s\right)$, as a function of $N$, on a logarithmic
scale. Only a sparse set of values of $N$ are sampled. The figure
on the left, which would confirm or disprove the conjecture, is ambiguous.
Although the data, as illustrated, shows a precipitous drop past the
critical line $\Re s=\frac{1}{2}$, it is not yet obvious whether
there will be a bounce/oscillation back up towards $\Re s=\frac{1}{2}$,
or whether the asymptotic behavior aims at $\Re s\simeq0.45$. Unfortunately,
computational costs are large; each additional point in these graphs
requires in excess of a week of compute time on current-era computers.

\lyxline{\normalsize}

\end{figure}

\section{Conclusions}

Numerical evidence suggests that the permutation-inspired zeta \[
\zeta_{p,q}(s)=\lim_{N\to\infty}\zeta_{N;p,q}(s)\]
obeys a Riemann hypothesis, in that it's zeros lie in the critical
strip, and particularly on the $\Re s=1/2$ line. High-precision numerical
computations are difficult, and offer only ambiguous support for the
conjecture. Curiously, the numerical evaluation itself {}``accidentally''
introduces some zeta-like artifacts, these artifacts do not seem to
be sufficiently strong to explain the overall structure. Alternate
approaches to the problem are desired, but the author is not aware
of any particular theoretical framework that would be applicable,
and could be employed to provide deeper insight into the phenomenon.

\bibliographystyle{plain}
\bibliography{/home/linas/linas/fractal/paper/fractal}

\end{document}